\UseRawInputEncoding 
 
\documentclass[12pt,twoside]{article}

\usepackage{graphicx}
\usepackage{amsfonts}

\begin{document}

\pagestyle{plain}

\date{}

\title{\bf Knot Dynamics}
\author{Louis H. Kauffman,\\ Department of Mathematics,\\ 851 South Morgan Street, \\ University of Illinois at Chicago,\\
Chicago, Illinois 60607-7045,\\ $<$kauffman@uic.edu$>$\\}

\maketitle

\noindent {\bf Abstract}
The paper studies the dynamics of knots under self-repulsion produced by artificial charges placed on the curve in space representing the knot.\\

\section{Introduction}
It is the purpose of this paper to ask questions about the dynamics of knots under self-repulsion.
A knot is arranged in three dimensional space so that it is coated with a mathematical version of electrical charge so that it evolves to 
geometric positions that minimize its energy under constraints that maintain its total length. Will such evolutions result in the unknotting of 
loops that are topologically unknotted? What are the forms of extermal positions in such evolutions?\\

We will give a guided tour of computer experiments that can be performed in the quest for answers to these questions. In the course of our specific computer exploration we use
the knot theory of rational knots and rational tangles \cite{HKCT} to produce classes of unknots with complex initial configurations that we call {\it hard unknots} and corresponding
complex configurations that are topologically equivalent to simpler knots. We shall see that these hard unknots and complexified knots give examples that do not reduce in the experimental
space of the computer program. That is, we find unknotted configurations that will not reduce to simple circular forms under self-repulsion, and we find complex versions of knots that will not 
reduce to simpler forms under the self-repulsion. It is clear to us that the phenomena that we have discovered depend very little on the details of the computer program as long as it 
conforms to our general description of self-repulsion. Thus, we suggest on the basis of our experiments that sufficiently complex examples of hard unknots and sufficiently complex examples
of complexified knots will not reduce to global minimal energy states in self-repulsion environments. In the course of the paper we will make the character of these examples precise.
It is our challenge to other program environments to verify or disprove our assertions.\\

To make this challenge quite specific, here is a family of examples. Let $$U_{n} = N([[n+1]] - [[n]]) = [(n+1).n]$$ in the notation of Section 3 of this paper. $U_{n}$ is a diagram for an unknot.
The example $[11.10]$ is shown in Figure~\ref{Unknot} and Figure~\ref{Unknot1}. We conjecture that any algorithm based only on self-repulsion will be unable to unknot $U_{n}$ for some natural number value for $n.$
In this example, our version of KnotPlot is unable to unknot $[11.10].$\\

More generally, again using the notation of Section 3, we have an infinity of unknots in the form $N(T-S)$ where $T=(a_1, a_2, \cdots , a_{n-1}, a_{n})$ is a  rational tangle in continued fraction form, and $S = (a_1, a_2, \cdots , a_{n-1})$ is the {\it truncate} of T obtained by removing the last term of $T.$ These unknotted configuratons provide a wealth of inputs for knot dynamical experiments.
In Section 3 we explain how to create complex configurations of given knot types by the same method.\\

Having made this challenge, we point out that in \cite{RS} refined methods of gradient descent provide ways in which many examples that we do not reach with our methods can be unknotted.
And in \cite{SchareinBFACF} the BFACF algorithm \cite{vanR} is used to strong effect. The challenge is given to all of these more powerful algorithms as well.\\

A knot is represented by a (cyclic) sequence $p_1,p_2,\cdots, p_n$ of points in three dimensional space
$R^3.$ In computer graphics the points can be shown and edges can connect point $p_i $ with point $p_{i+1}$ and of course $p_n$ with $p_1 .$ Another computer
graphic option is to enclose the resulting piecewise linear curve with  tube so that the knot appears as a smooth tube with a chosen diameter. We make use of the computer program
KnotPlot \cite{KnotPlot} created by Robert Scharein, and we thank him for numerous conversations about the use of this program, and for modifying the program to facilitate the use of
tangles and rational knots. We shall refer to KnotPlot when speaking of experiments
and examples produced in the course of the paper.\\

A knot undergoes dynamics of self-repulsion in the computer program  by a recursion involving forces of order $r^{-5}$ ( $r$ is the distance between a chosen pair of points on the knot) applied to all pairs of non-adjacent points in the sequence of points that define the knot (see above).
The recursion is controlled so that the knot-type of the space curve does not change from step to step. We can calcluate the Simon Energy \cite{Simon}. The Simon Energy is the sum of $1/|p-q|$ over all pairs of distinct points in the list of points in three dimensional space that describes the knot.
The KnotPlot program minimizes this energy by self-repulsion forces of order $1/|p-q|^5 .$ The Simon Energy corresponds directly to forces of the form $1/|p-q|^2 .$ \\

One can examine many aspects of the knot dynamics. Knots and links evolve into special geometric configurations that are not predictable beforehand.
One often finds that unknotted curves will unfold and become simplified into curves that are nearly flat circles. In some cases the force field cannot undo unknots.
The final  configuration can depend upon the starting configuration and upon degrees of perturbation that are allowed in the process of evolution. 
The program we use allows, in addition to the
repulsion forces, an attraction of nearby points via spring-forces (using an analog of Hooke's Law) \cite{SchareinThesis}. 
The balance of the spring and repulsion forces can be used to facilitate experiments.
We will give examples of all these phenomena.\\

\section{Moving Down Energy Gradients}
We begin with a relatively flat representation of the trefoil knot, a torus knot of type $(2,3).$ Generally we consider torus knots of type $(a,b)$ where $a$ and $b$ are relative prime positive
integers. Such a knot winds around a standardly embedded torus in three-space $a$ times around the meridian direction and $b$ times around the longitude direction. It is a fact that a torus knot
of type $(a,b)$ is topologically equivalent to a torus knot of type $(b,a).$ We shall write $T(a,b)$ for a torus knot of type $(a,b)$ and so write that $T(a,b)$ and $T(b,a)$ are topologically equivalent. We shall see this equivalence emerge in the dynamics.\\

\begin{figure}
     \begin{center}
     \begin{tabular}{c}
     \includegraphics[width=7cm]{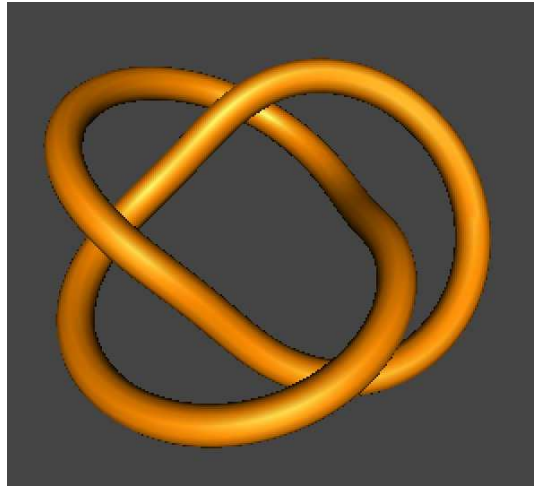}
     \end{tabular}
     \caption{\bf Initial Trefoil}
     \label{T1}
\end{center}
\end{figure}

\begin{figure}
     \begin{center}
     \begin{tabular}{c}
     \includegraphics[width=7cm]{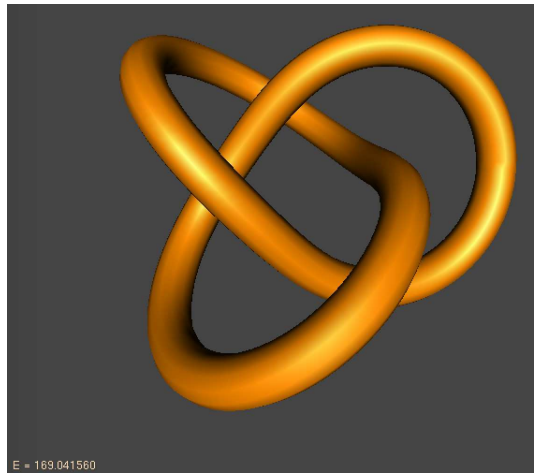}
     \end{tabular}
     \caption{\bf Self Repelled Trefoil}
     \label{T2}
\end{center}
\end{figure}

\begin{figure}
     \begin{center}
     \begin{tabular}{c}
     \includegraphics[width=7cm]{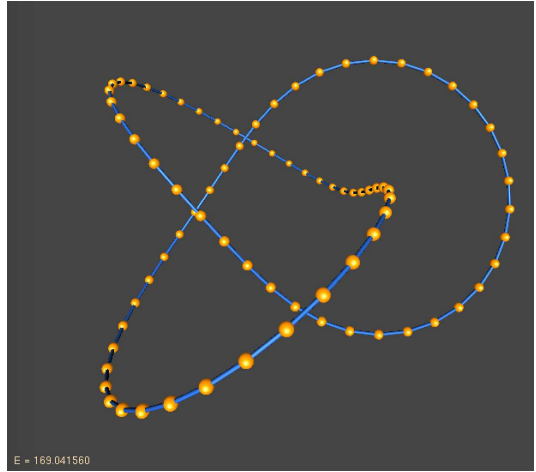}
     \end{tabular}
     \caption{\bf Beaded Trefoil}
     \label{T3}
\end{center}
\end{figure}

In Figure~\ref{T1} we see the initial representation of $T(2,3).$ In Figure~\ref{T2} we see the result of the trefoil after it has reached self-repulsion stability. At this point the trefoil configuration has
calculated Simons energy $E = 169.041560.$  The Simons energy is computed by using the fact that the knot we see is represented by a collection of points in three dimensional Euclidean space as shown in Figure~\ref{T3}. The Simons energy is the sum of $1/r(p,q)$ ranging over all pairs of points $(p,q)$ in this configuration where $r(p,q)$ is the Euclidean distance between the points
$p$ and $q.$ The Simon energy corresponds classically to the force field with exponent $2,$ but it is convenient to track this energy even when using other force functions.\\

\begin{figure}
     \begin{center}
     \begin{tabular}{c}
     \includegraphics[width=7cm]{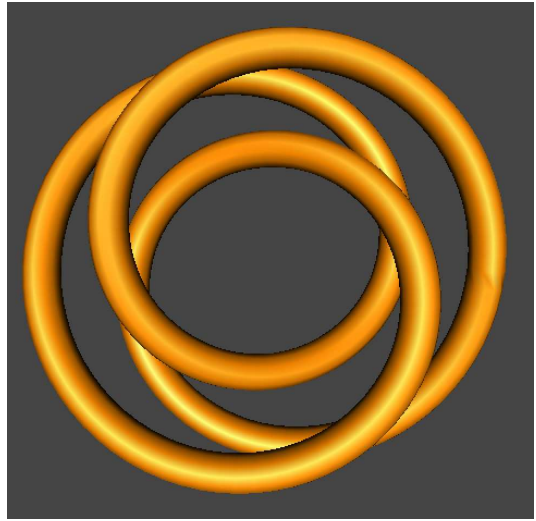}
     \end{tabular}
     \caption{\bf (3,2) Torus Knot}
     \label{TT1}
\end{center}
\end{figure}

\begin{figure}
     \begin{center}
     \begin{tabular}{c}
     \includegraphics[width=7cm]{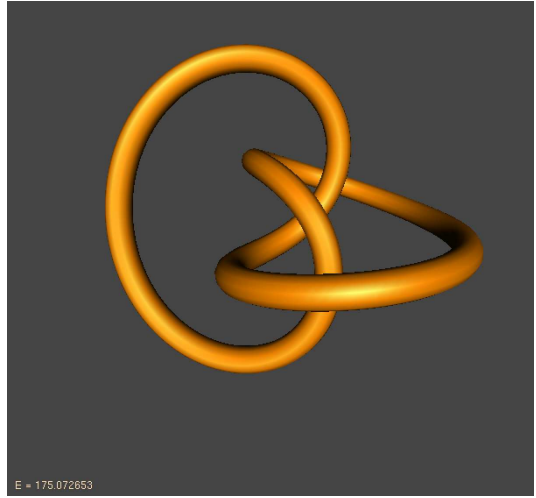}
     \end{tabular}
     \caption{\bf Self Repelled (3,2) Torus Knot}
     \label{TT2}
\end{center}
\end{figure}

\begin{figure}
     \begin{center}
     \begin{tabular}{c}
     \includegraphics[width=7cm]{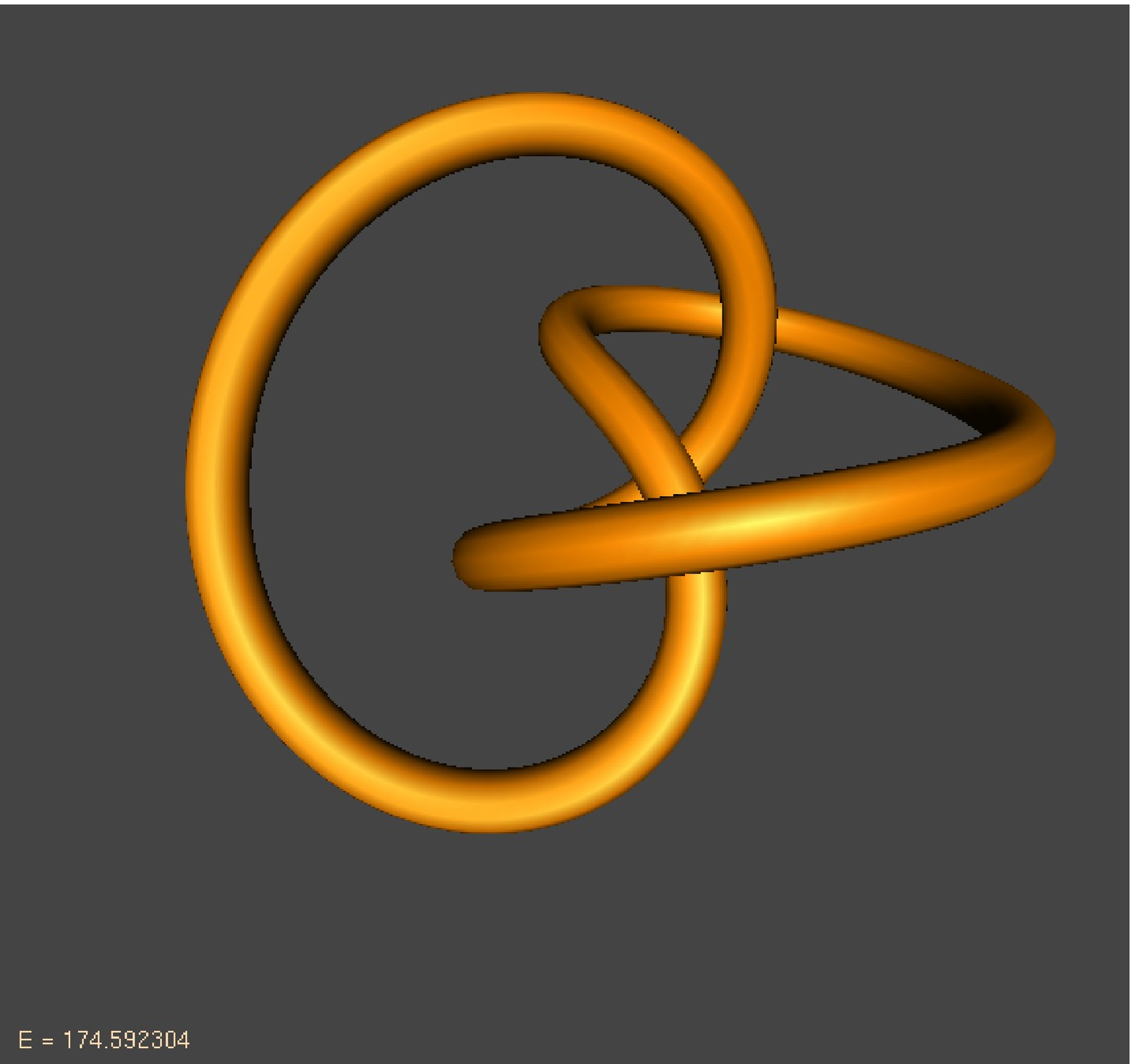}
     \end{tabular}
     \caption{\bf Spring Dynamics Destabilizes the (3,2).}
     \label{TT3}
\end{center}
\end{figure}

\begin{figure}
     \begin{center}
     \begin{tabular}{c}
     \includegraphics[width=7cm]{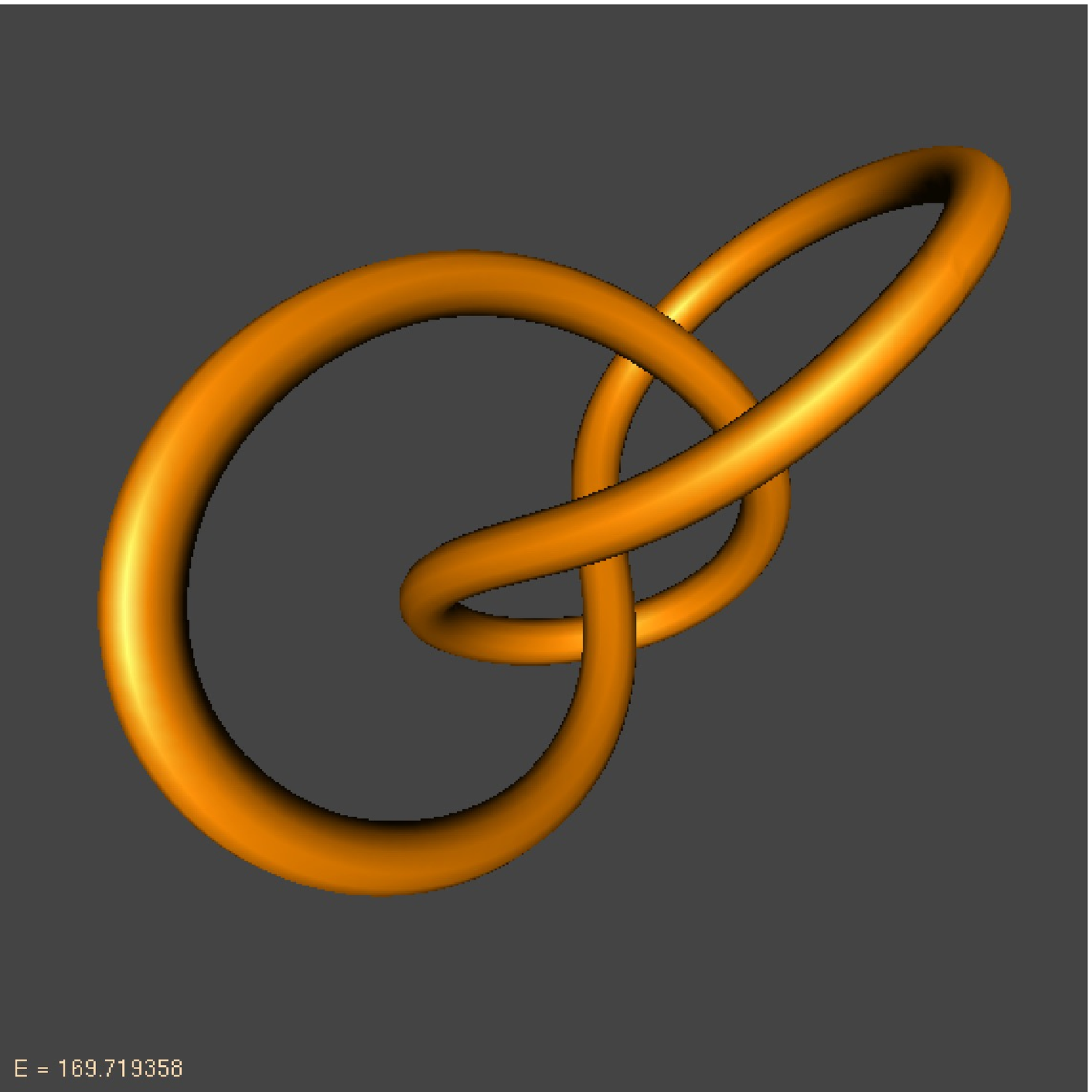}
     \end{tabular}
     \caption{\bf Spring Dynamics Destabilizes the (3,2).}
     \label{TT4}
\end{center}
\end{figure}

\begin{figure}
     \begin{center}
     \begin{tabular}{c}
     \includegraphics[width=7cm]{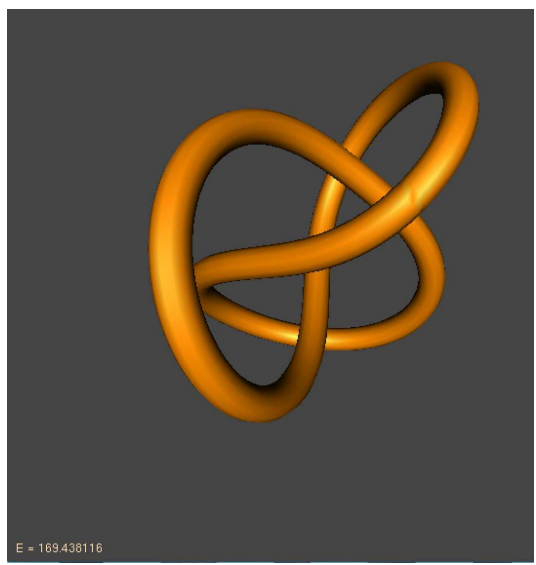}
     \end{tabular}
     \caption{\bf Toward a New Minimum.}
     \label{TT5}
\end{center}
\end{figure}

\begin{figure}
     \begin{center}
     \begin{tabular}{c}
     \includegraphics[width=7cm]{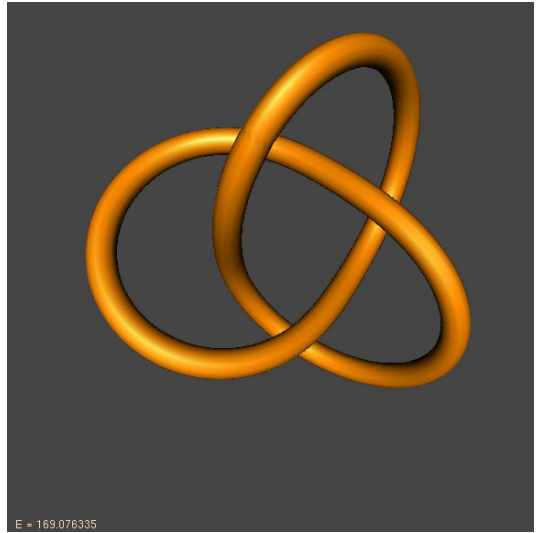}
     \end{tabular}
     \caption{\bf Final Minimum Equivalent to the (2,3) Trefoil Minimum.}
     \label{TT6}
\end{center}
\end{figure}

In Figure~\ref{TT1} we see the initial state for a torus knot of type $T(3,2).$ Topologically, a torus knot of type $T(a,b)$ is equivalent to one of type $T(b,a),$  but geometrically this is quite a different starting position. In one case the curve goes around the torus longitudinally twice and meridianally three times. This is reversed to three and two in the second case. We now examine how the repulsion dynamics worked on the $T(3,2).$ In Figure~\ref{TT2} we see the self-repelled $(3,2)$ torus knot in what appears to be a stable position. Experimentally this configuration is 
stable when using the undamped program in KnotPlot. Undamped means that the springs that are available between successive nodes in the model are not operation (in the sense that the spring constant is high and the springs are hence fully contracted). The energy level of this configuration is $174.59304$, considerably higher than the minimum we found for the $(2,3)$ knot. Now we have the opportunity to perturb this configuration by using the program in its undamped form. Then the springs are free to contract and expand, letting the energy move into kinetic energy of motion for the configuration. In Figures~\ref{TT3} ~\ref{TT4} ~\ref{TT5} ~\ref{TT6} we see that indeed the undamped evolution results in an oscillation that swings arcs of the knot around and lets it fall into the lower energy configuration of the $(2,3)$ trefoil knot.\\

\begin{figure}
     \begin{center}
     \begin{tabular}{c}
     \includegraphics[width=7cm]{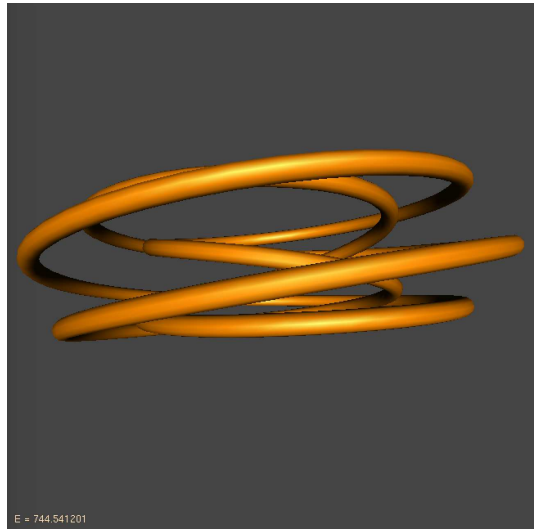}
     \end{tabular}
     \caption{\bf  Initial condition of $(5,2)$ torus knot.}
     \label{TF1}
\end{center}
\end{figure}

\begin{figure}
     \begin{center}
     \begin{tabular}{c}
     \includegraphics[width=7cm]{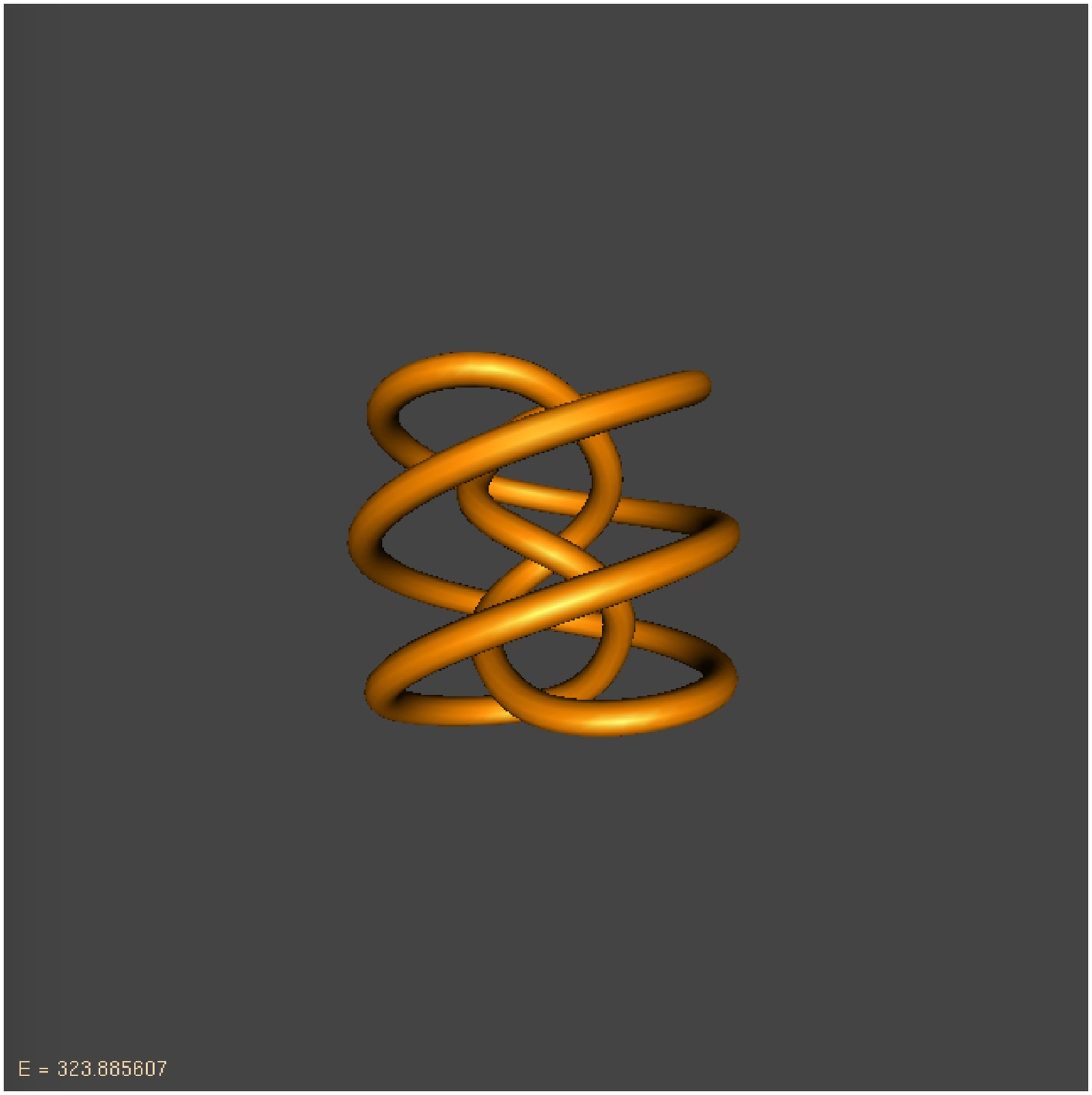}
     \end{tabular}
     \caption{\bf Descending $(5,2)$ torus knot.}
     \label{TF2}
\end{center}
\end{figure}

\begin{figure}
     \begin{center}
     \begin{tabular}{c}
     \includegraphics[width=7cm]{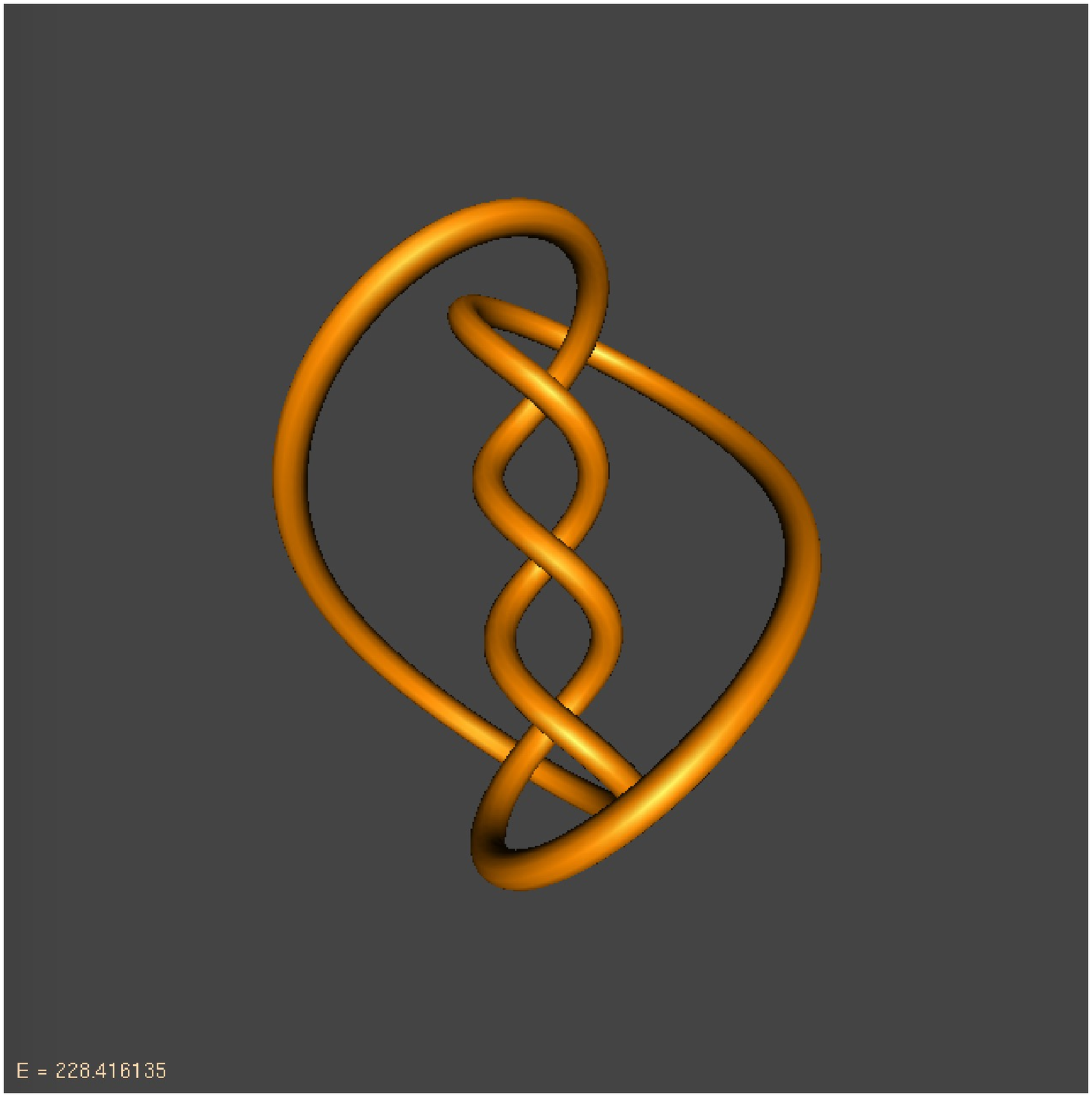}
     \end{tabular}
     \caption{\bf Descending $(5,2)$ torus knot.}
     \label{TF3}
\end{center}
\end{figure}

\begin{figure}
     \begin{center}
     \begin{tabular}{c}
     \includegraphics[width=7cm]{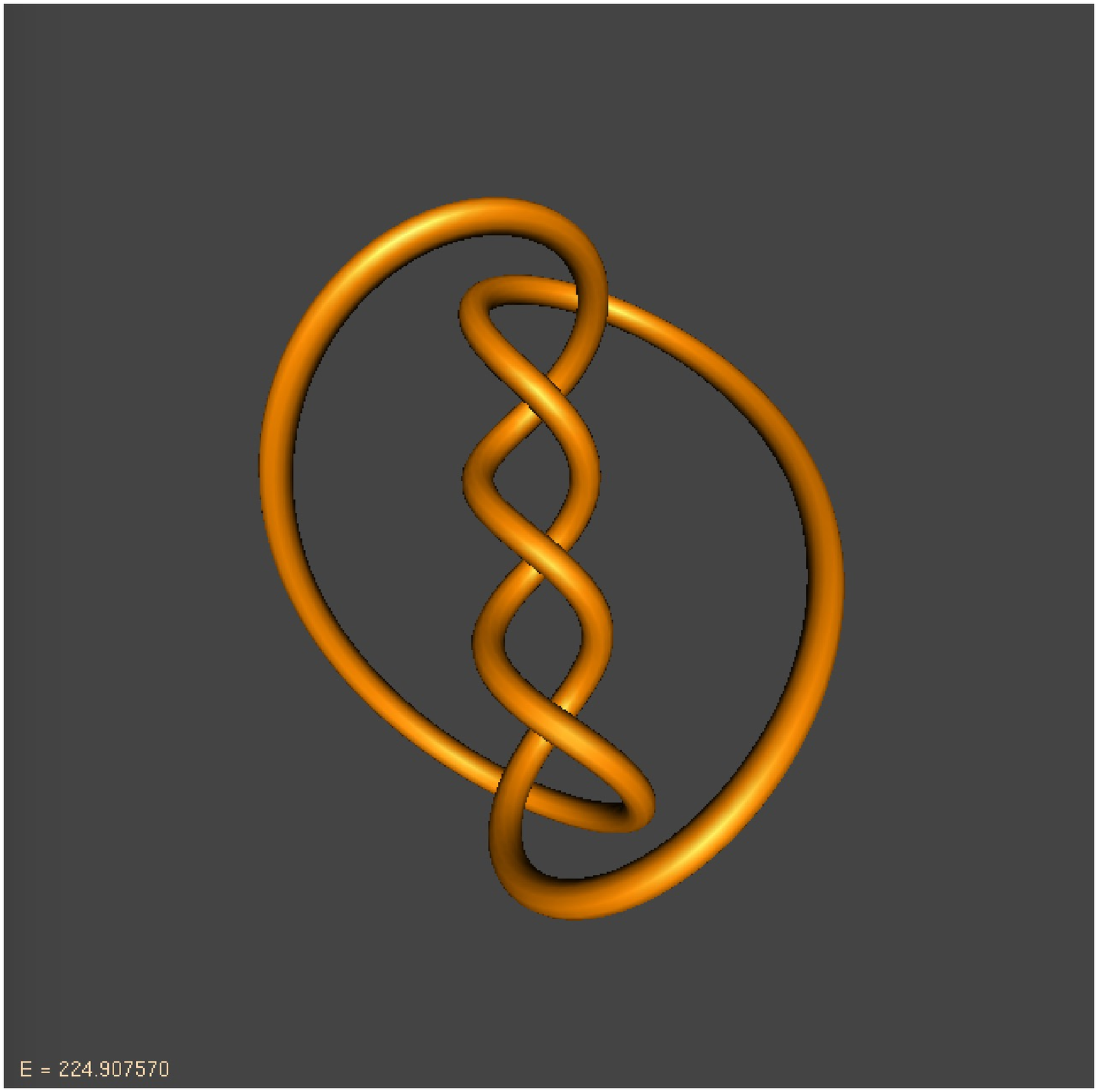}
     \end{tabular}
     \caption{\bf Descending $(5,2)$ torus knot.}
     \label{TF4}
\end{center}
\end{figure}

\begin{figure}
     \begin{center}
     \begin{tabular}{c}
     \includegraphics[width=7cm]{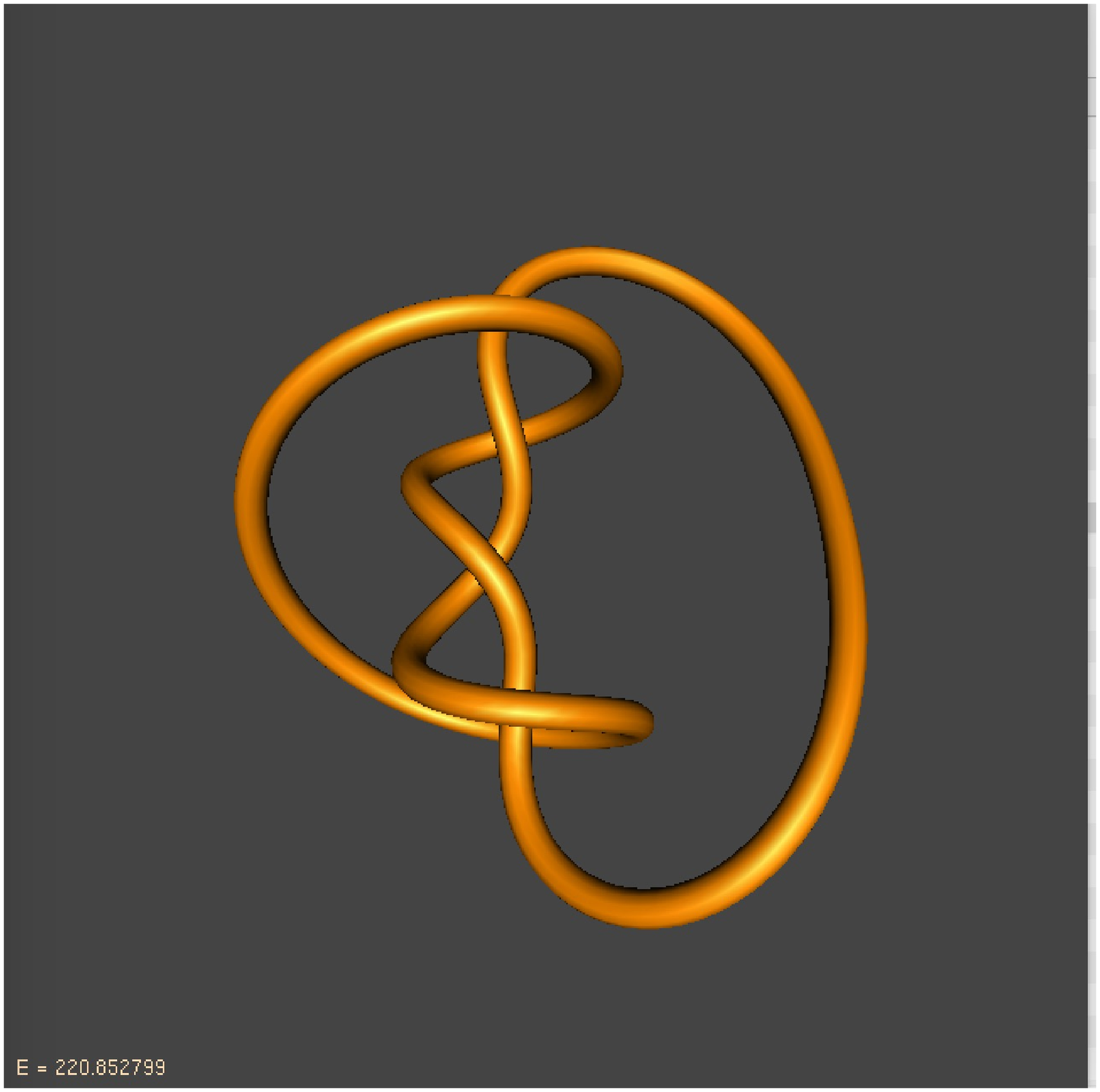}
     \end{tabular}
     \caption{\bf Descending $(5,2)$ torus knot.}
     \label{TF5}
\end{center}
\end{figure}

\begin{figure}
     \begin{center}
     \begin{tabular}{c}
     \includegraphics[width=7cm]{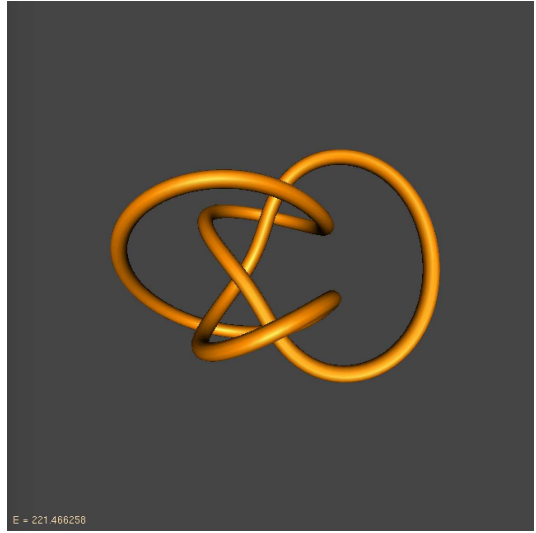}
     \end{tabular}
     \caption{\bf Final Descent of  $(5,2)$ torus knot.}
     \label{TF6}
\end{center}
\end{figure}

\begin{figure}
     \begin{center}
     \begin{tabular}{c}
     \includegraphics[width=7cm]{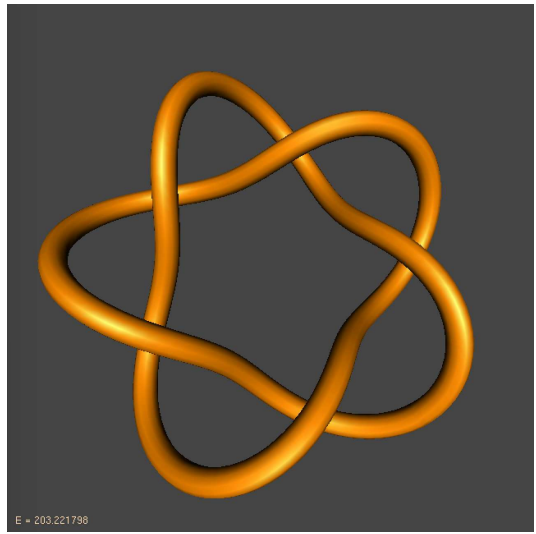}
     \end{tabular}
     \caption{\bf Minimum for the (2,5) torus knot.}
     \label{TF7}
\end{center}
\end{figure}

Figures~\ref{TF1} ~\ref{TF2} ~\ref{TF3} ~\ref{TF4}~\ref{TF5} ~\ref{TF6}  illustrate the descent of a $(5,2)$ torus knot. The resulting energy minimum is $221.466258.$ Note that the internal twist
in the knot collapes in the course of this evolution. The same collapse will be observed for $(n,2)$ torus knots with $n \ge 5.$
If one uses the undamped descent for a long time, then for the  $(5,2)$ the resulting configuration does not change when the undamped evolution is made available.
If one runs the damped evolution for a short time, then the result will be similar to our evolution with the trefoil knot and an exchange of kinetic energy takes the knot to the same configuration as
the descent of the $(2,5)$ torus knot. We illustrate the final descent configuration for a $(2,5)$ torus knot in Figure~\ref{TF7}.\\

\begin{figure}
     \begin{center}
     \begin{tabular}{c}
     \includegraphics[width=7cm]{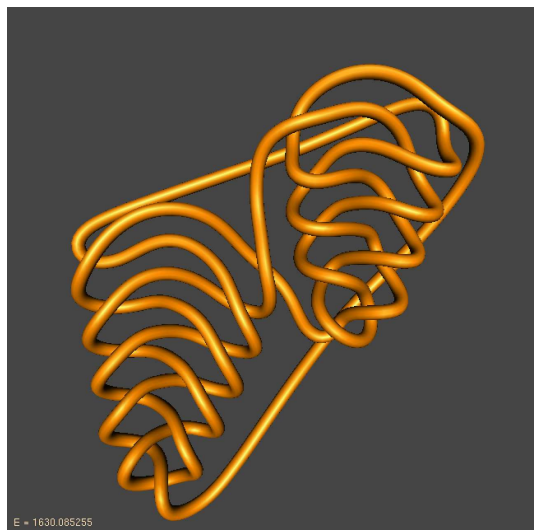}
     \end{tabular}
     \caption{\bf Starting Position for the unknot $[11.10]$.}
     \label{Unknot}
\end{center}
\end{figure}

\begin{figure}
     \begin{center}
     \begin{tabular}{c}
     \includegraphics[width=7cm]{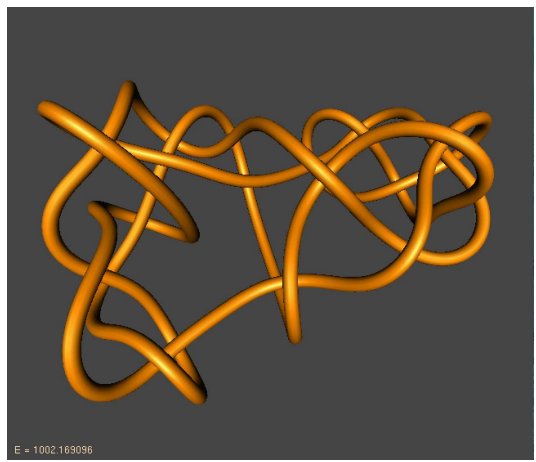}
     \end{tabular}
     \caption{\bf End Position for the unknot $[11.10]$ in Local Energy Minimization - A high energy version of the unknot.}
     \label{Unknot1}
\end{center}
\end{figure}

In Figure~\ref{Unknot} we illustrate an initial position for an unknot of type $[11.10]$ where we will explain this notation below. The Figure~\ref{Unknot1} shows the result of the undamped evolution of this knot. It has not become unknotted and no application of the undamped program will unknot it at any point in the process. This is an eperimental example of an unknnotted curve that fails to become unknotted in the force field of the Knotplot program. We conjecture that this phenomenon is not specific to this program but rather that some unknots of the form 
$[(n+1).n]$ will not be able to evolve to simple circles in any forcefield of this kind. The reader is advised to experiment with these examples. We explain the corresponding tangle theory in the next section.\\

\begin{figure}
     \begin{center}
     \begin{tabular}{c}
     \includegraphics[width=7cm]{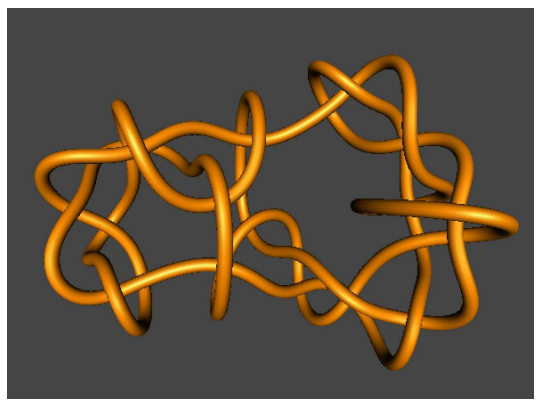}
     \end{tabular}
     \caption{\bf End Position for the knot $[15.10]$ in Local Energy Minimization - A high energy version of the figure eight knot.}
     \label{KnotMin}
\end{center}
\end{figure}

The knot shown in Figure~\ref{KnotMin}  is obtained by the same techniques as the previous two figures. It is of type $[15.10]$ and is equivalent to a figure eight knot as we will explain in the next section.It will not reduce any further in the KnotPlot environment either via undamped or damped evolution. In fact, as we will see in then next section, we have an infinite collection of figure eight knots in the forms $[(n+5).n]$ and for $n$ greater than $10$ we have the empirical result that they all have local energy minima , increasing with the value of $n$ in the KnotPlot environment.
We conjecture that this phenomenon is qualitatively the same in any force field environment for knot repulsion.\\

\clearpage

\section{Hard Knots and Collapsing Tangles}

\begin{figure}
     \begin{center}
     \begin{tabular}{c}
     \includegraphics[width=9cm]{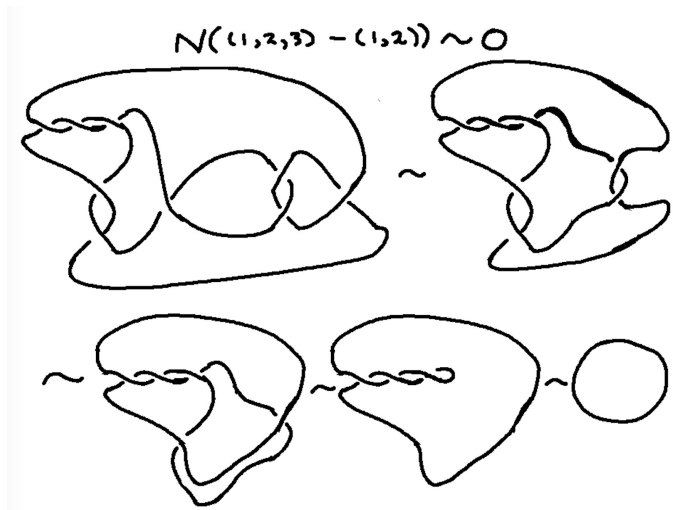}
     \end{tabular}
     \caption{\bf $N((1,2,3) - (1,2))$ is Unknotted.}
     \label{UKK3}
\end{center}
\end{figure}

Before we begin the details of this section, the reader should view Figure~\ref{UKK3}. There is illustrated a diagram and a sequence of steps that show that the diagram is
unknotted. This is the type of unknot that we have discussed in the previous section. In order to systematically construct configurations of this type, that are complex and yet unknotted, we
develop an algebraic language of tangles in the present section of the paper. For complete information about all aspects of this discussion we recommend that the reader consult
the paper \cite{HKCT}.\\

In this section we give background in terms of a theory of tangles for constructing the unknots and knots in collapsible states.
We begin with the definition of a tangle and the notions of combination of tangles (addition, star product, numerator and denominator) and how this leads to the tangle fraction, the notion of 
rational tangles and rational knots and how these constructions are related to computer experiments.The notion of a tangle was introduced in 1967 by
Conway \cite{Conway} in his work on enumerating and classifying knots and links. 
\\

\begin{figure}
     \begin{center}
     \begin{tabular}{c}
     \includegraphics[width=10cm]{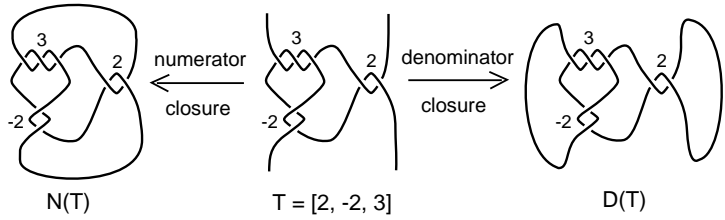}
     \end{tabular}
     \caption{\bf Tangles - Numerators and Denominators}
     \label{UK9}
\end{center}
\end{figure}

\begin{figure}
     \begin{center}
     \begin{tabular}{c}
     \includegraphics[width=10cm]{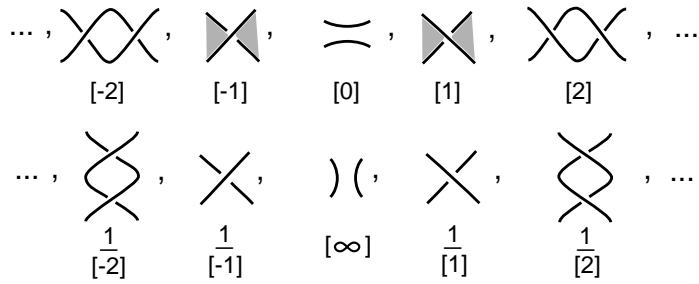}
     \end{tabular}
     \caption{\bf Elementary Rational Tangles}
     \label{UK10}
\end{center}
\end{figure}

A {\it $2$-tangle} is a proper embedding of two unoriented arcs and a finite number of circles in a $3$-ball
$B^3,$ so that the four endpoints lie in the boundary of $B^3$. A {\it tangle diagram} is a regular projection
of the tangle on an equatorial disc of $B^{3}.$   A {\it rational tangle} is a special case of a $2$-tangle obtained by applying consecutive twists
on neighbouring endpoints of two trivial arcs. Such a pair of arcs comprise the $[0]$ or $[\infty]$ tangles 
as shown in Figure~\ref{UK10}, depending on their position in the plane. We shall say that the rational tangle is in {\it
twist form} when it is obtained by such successive twists of adjacent strands. For example see Figure~\ref{UK9} where we show a tangle in the center of the figure that is obtained by 
first starting with a horizontal 3-twist ($[3]$ in Figure~\ref{UK10} and then twisting two vertical strands negatively twice and then twisting two horizontal strands positively twice.
Tangles obtained recursively from the untwisted tangles $[0]$ and $[\infty]$ of Figure~\ref{UK10} are called {\it rational tangles}.\\

\begin{figure}
     \begin{center}
     \begin{tabular}{c}
     \includegraphics[width=7cm]{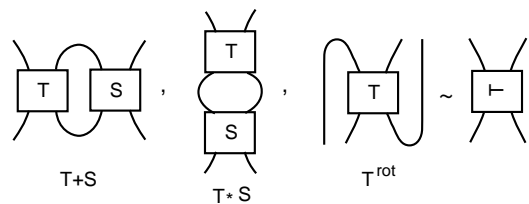}
     \end{tabular}
     \caption{\bf Addition, Star Product and Rotation of Tangles}
     \label{UK11}
\end{center}
\end{figure}

Note that in Figure~\ref{UK10} we have indicated horizontal twists by bracketed integers $[n]$ and vertical twists by reciprocal integers $1/[n],$ including the convention that 
$1/[0] = [\infty].$ In Figure~\ref{UK11} we show tangle operations $T+ S$, $T \star S$, and $T^{rot}.$ These are called {\it addition, star product} and {\it rotation} of tangles. The rotation operation has the effect of turning the tangle  by ninety degrees in the plane. Note that twisting horizontal strand $n$ times is the same as forming $T + [n]$ and twisting vertical strands $n$ times is the same as forming $T \star [n]$. Thus we can describe the tangle in Figure~\ref{UK9} as ($[3] \star [-2])+[2]$. \\

It turns out that one can associate to any tangle $T$ a fraction $F(T)$ with values in the rational numbers plus the formal value $\infty,$ and that this fraction satisfies the following rules with
respect to addition, star product and rotation. Here $T^{\star}$ denote the mirror image tangle, obtained by switching all the crossings of $T.$
Once we assume that $F([0]) = 0, F([1]) = 1,$ these rules determine the fractions of rational tangles.\\

\begin{enumerate}
\item $F(S+T) = F(S) + F(T)$
\item $F(T \star S) = \frac{1}{\frac{1}{F(T)}+ \frac{1}{F(S)}}$
\item $F(T^{rot}) = -1/F(T)$
\item $F(T^{\star}) = - F(T)$
\end{enumerate}

$~$\\
\noindent The second property of the fraction for the star product is a conseqence of first property for addition, using the third property that tells us that the fraction of a rotate is the negative reciprocal of 
the fraction of the original tangle. Note that it follows from these rules that $F([n]) = n$ and that $F(1/[n]) = 1/n$ and that 
$$F(T \star (1/[n])) = 1/(n + 1/F(T)).$$
This is the formula for the fraction of a tangle with an extra vertical twist $1/[n]$ at the bottom.
From this formula we see that the fraction of a rational tangle will be represented by  a continued fraction. For example,
$$F([a] + ([c] \star (1/[b]))) =  a + 1/( b + 1/c).$$
For this reason it is useful to denote rational tangles by just writing a corresponding continued fraction. 

\begin{figure}
     \begin{center}
     \begin{tabular}{c}
     \includegraphics[width=10cm]{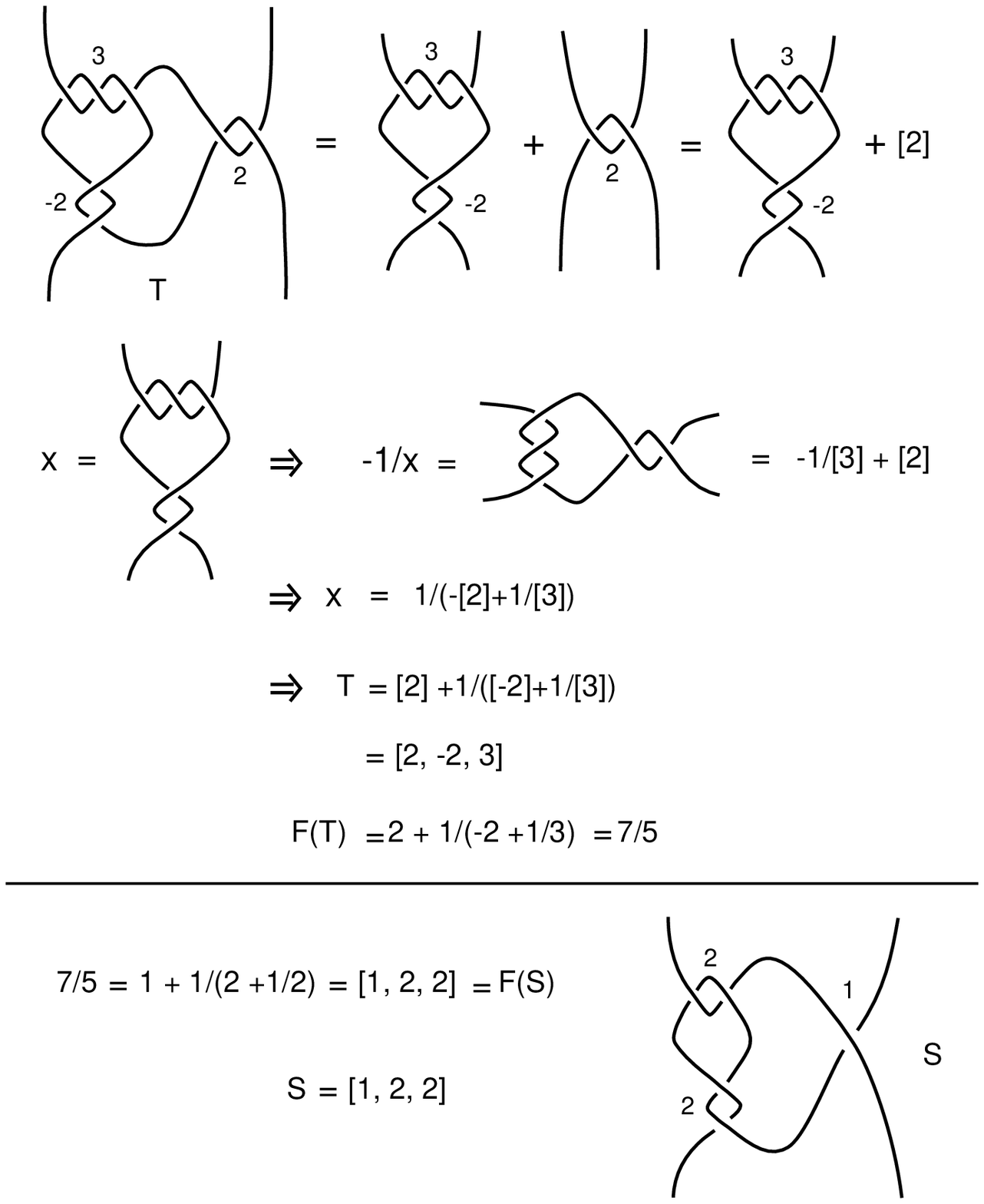}
     \end{tabular}
     \caption{\bf Finding the Fraction of a Rational Tangle}
     \label{UK12}
\end{center}
\end{figure}

\begin{figure}
     \begin{center}
     \begin{tabular}{c}
     \includegraphics[width=9cm]{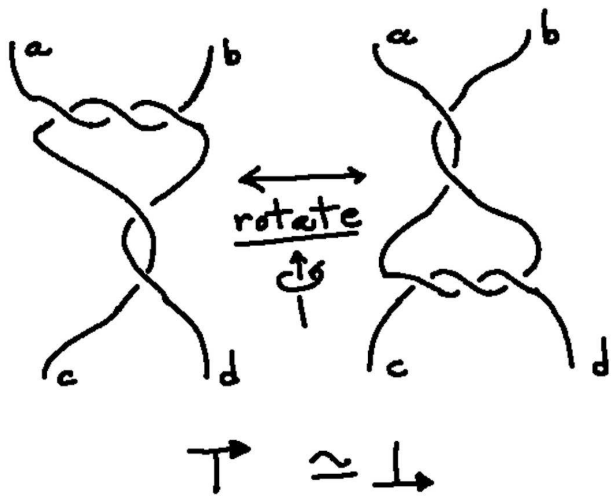}
     \end{tabular}
     \caption{\bf Rotational Isotopy}
     \label{UKK1}
\end{center}
\end{figure}

\begin{figure}
     \begin{center}
     \begin{tabular}{c}
     \includegraphics[width=9cm]{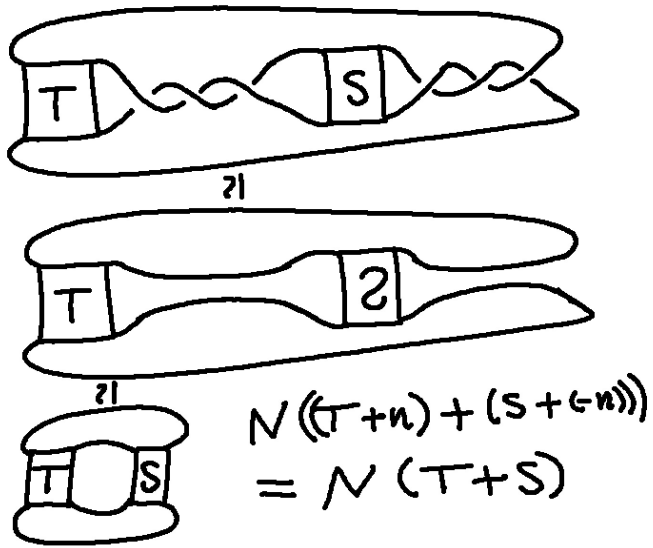}
     \end{tabular}
     \caption{\bf Simplifying a Tangle Sum}
     \label{UKK2}
\end{center}
\end{figure}

\begin{figure}
     \begin{center}
     \begin{tabular}{c}
     \includegraphics[width=9cm]{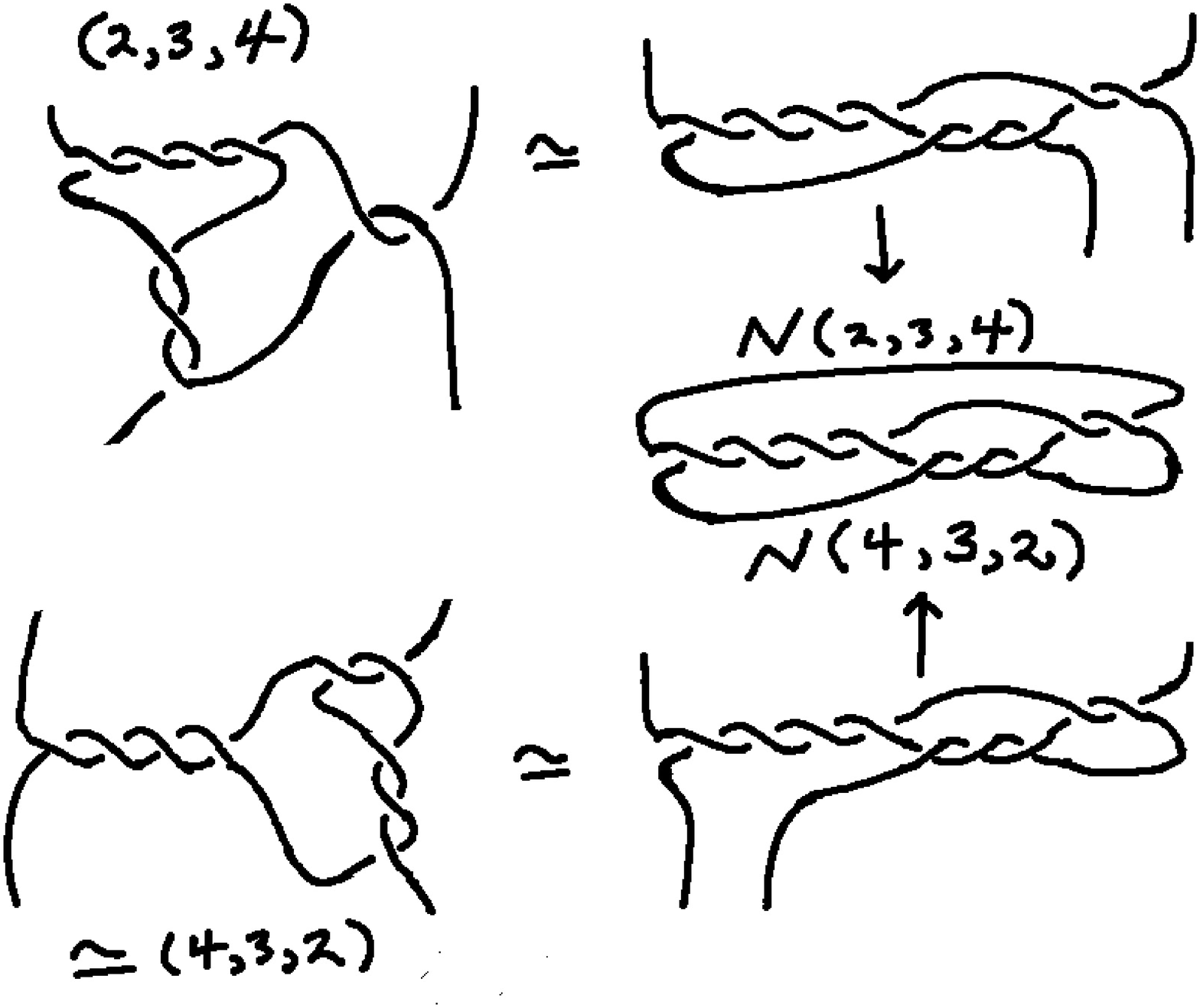}
     \end{tabular}
     \caption{\bf $N(a,b,c) = N(c,b,a)$}
     \label{UKK4}
\end{center}
\end{figure}

\begin{figure}
     \begin{center}
     \begin{tabular}{c}
     \includegraphics[width=10cm]{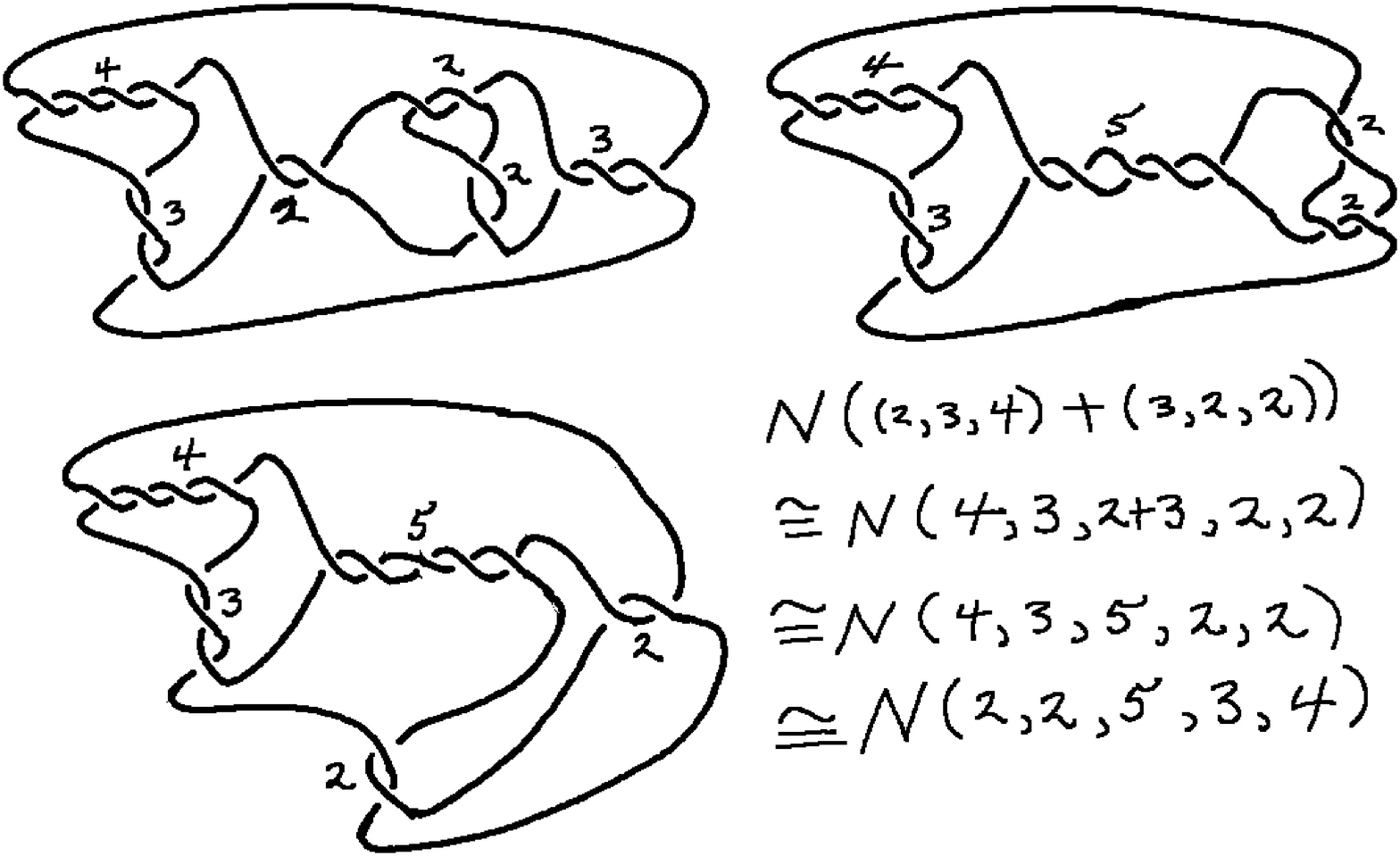}
     \end{tabular}
     \caption{\bf $N((a,b,c) + (d,e,f)) = N(c,b,a+d,e,f)$}
     \label{UKK5}
\end{center}
\end{figure}

\begin{figure}
     \begin{center}
     \begin{tabular}{c}
     \includegraphics[width=10cm]{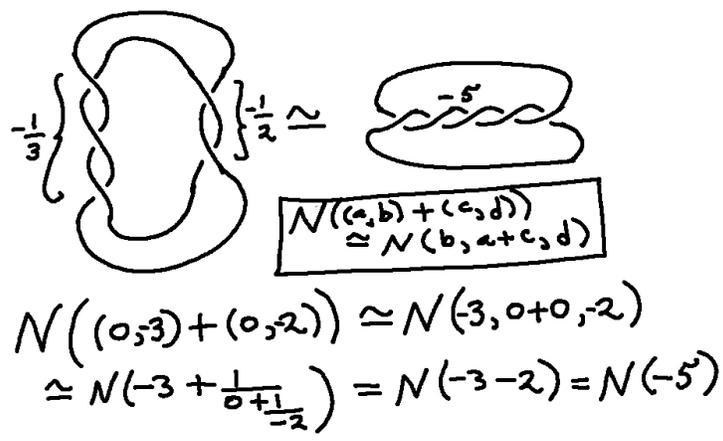}
     \end{tabular}
     \caption{\bf Second example of the numerator of a sum of tangles.}
     \label{UKK6}
\end{center}
\end{figure}

 \clearpage

\begin{figure}
     \begin{center}
     \begin{tabular}{c}
     \includegraphics[width=10cm]{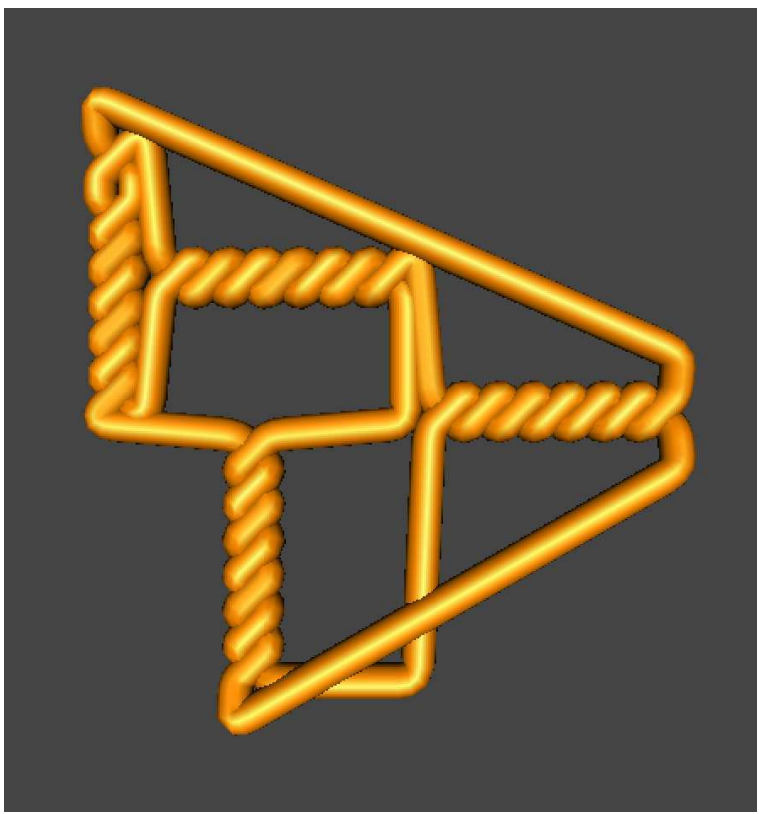}
     \end{tabular}
     \caption{\bf  $N(7,7,7,7)$ -- Initial State.}
     \label{UKK7}
\end{center}
\end{figure}

\begin{figure}
     \begin{center}
     \begin{tabular}{c}
     \includegraphics[width=10cm]{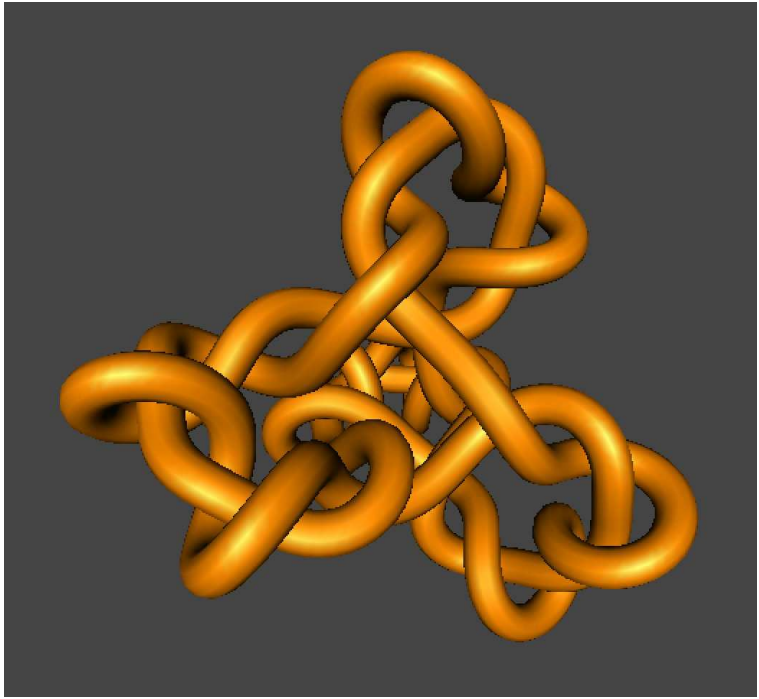}
     \end{tabular}
     \caption{\bf  $N(7,7,7,7)$ -- Final State.}
     \label{UKK8}
\end{center}
\end{figure}

  \begin{figure}
     \begin{center}
     \begin{tabular}{c}
     \includegraphics[width=10cm]{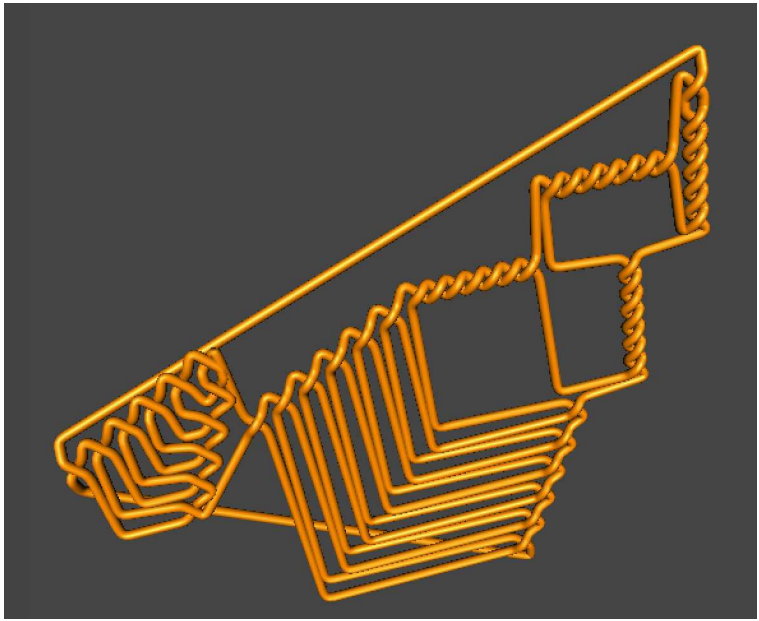}
     \end{tabular}
     \caption{\bf  Last Knot -- Initial State.}
     \label{UKK9}
\end{center}
\end{figure}

\begin{figure}
     \begin{center}
     \begin{tabular}{c}
     \includegraphics[width=10cm]{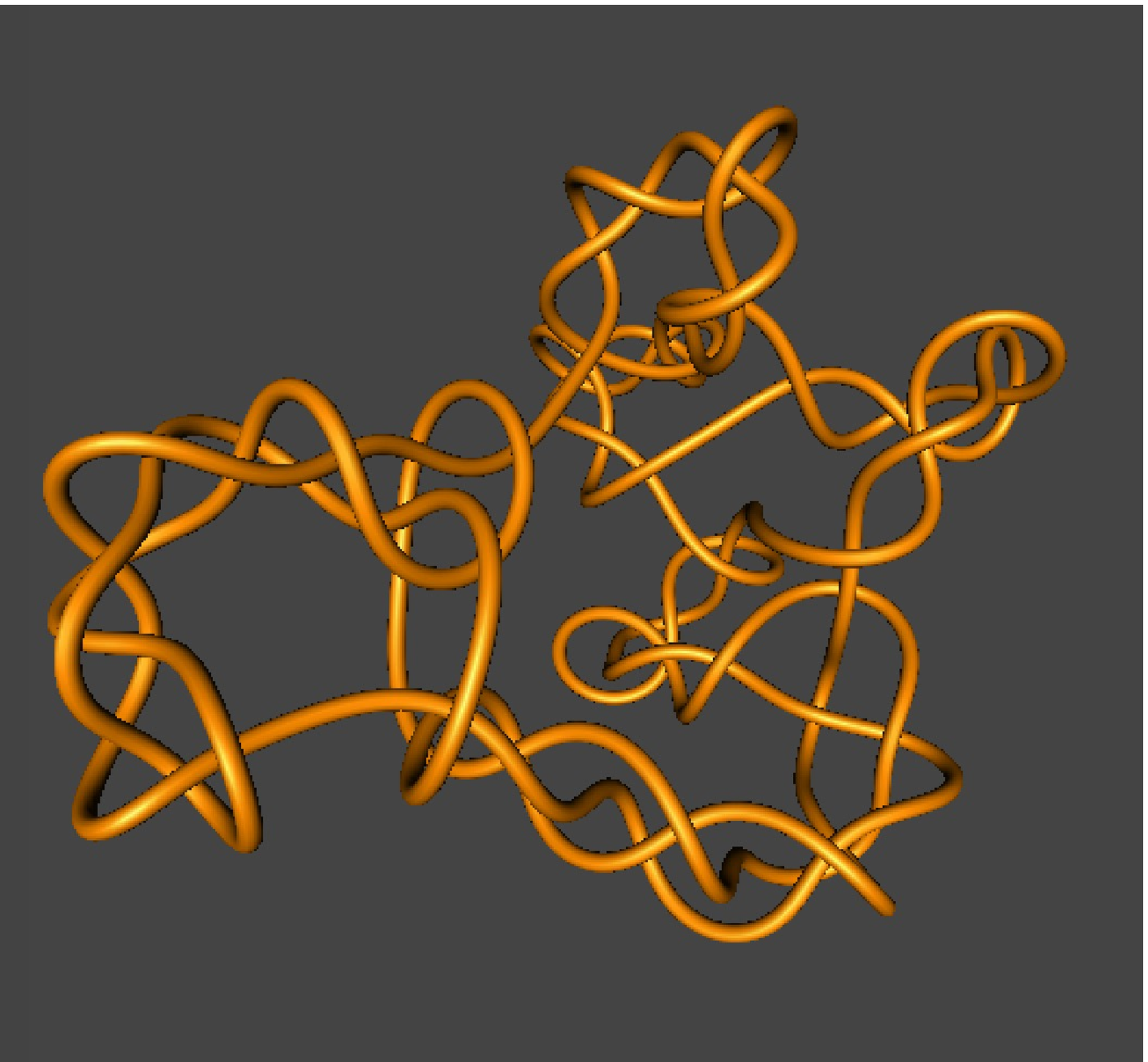}
     \end{tabular}
     \caption{\bf  Last Knot -- Final State.}
     \label{UKK10}
\end{center}
\end{figure}

We will illustrate this with an example.
In Figure~\ref{UK12} we illustrate these properties by starting with the tangle $T$ from Figure~\ref{UK9} and noting that $T$ is a sum of the form 
$T = S + [2]$ where $S = [3] \star [-2].$ We could directly use the rules above to compute as follows.
$$F(S) = F([3] \star F(1/[-2]) =  1/(1/F([3]) + F([-2])) $$
$$= 1/(1/3 + -2)) = 3/(1-6) = -3/5.$$
$$F(T) = F(S + [2]) = F([S]) + 2 = -3/5 + 2 = 7/5.$$
In the Figure~\ref{UK12} we denote $x = F(S)$ and use property $3.$ above to rotate the tangle and then use property $1.$ to complete the calculation of $F(S)$ without using property $2.$
In the rest of the figure we note that $7/5$ has a positive continued fraction expansion as $7/5 = 1 + 2/7 = 1 + 1/(7/2) = 1 + 1/(2 + 1/2) = [1,2,2].$ The last symbol, $[1,2,2]$ is the symbol for this continued fraction. Thus $$[a,b,c] = a + 1/(b + 1/c)$$ and more generally
$$[a_1, a_2, a_3 , a_4 , \cdots , a_n] = a_1 + 1/(a_2 + 1/(a_3 + 1/(a_4 + \cdots + 1/a_n)))).$$\\

Thus we shall from now on use the continued fraction written numerically to indicate the corresponding rational tangle. Thus $(1,2,2)$ will be our notation for the tangle diagram at the bottom of
Figure~\ref{UK12}. Note that topologically this tangle is equivalent to $(1,3,-2).$ These symbols represent the same numerical fraction and they each represent specific diagrams for rational tangles. The Figure~\ref{UK12} gives the reader enough information to see that these two tangles are topologically equivalent. That is, the reader can take as an exercise to deform the tangle
$(1,3,-2)$ to the tangle $(1,2,2)$  by an isotopy that does not move any arcs over the endpoints of the tangle, while leaving the endpoints fixed. \\

John Horton Conway proved that {\it two rational tangles are topolgically equivalent if and only if they have the same fraction.} See \cite{Conway,HKCT,KL1,KL2} for statements and proofs of the Conway Theorem. The key to our approach to proving the Conway Theorem, in these cited papers, is the insight that any continued fraction tangle with a mixture of signs such as  $(1,3,-2)$  
is not an alternating weave, and can be simplified by standard isotopies (swinging arcs that go over or under consecutive crossings)  until it is alternating. The rational form of the tangle is 
retained under such isotopies, and the final alternating form corresponds to an entirely positive (or entirely negative)  sequence  such as  $(1,2,2).$ The only positive or only negative  sequence is uniquely obtained from the fraction itself
by Euclid's algorithm. We give no more details here, but refer to our papers.\\

Note that $a + 1/(b + 1/c) = (a,b,c).$ Thus $(0,b,c) = 1/(b + 1/c) = 1/(b,c).$ These sequential symbols are useful for handling rational tangles and continued fractions. Note also that when we 
draw $(a,b,c)$ the twists occur from right to left. Thus in Figure~\ref{UK12} we illustrate $(2,-2, 3)$ and the $2$ occurs on the right in the diagram.\\

In Figure~\ref{UKK1} we illustrate a useful fact about about tangles. If $T$ is a tangle and we apply rotation to it twice it will be rotated by 180 degrees. It turns out that for rational tangles, one can pick them up, rotate them by 180 degrees and they are topologically unchanged. The figure illustrates one case of the internal topological rotations that accomplish this fact.\\

The next operation we need is the {\it numerator construction} $N(T)$ that converts a tangle to a knot of link by connect the top strands to one another and the bottom strands to one another.
See Figure~\ref{UK9} for an illustration of the numerator operation and its so-called corresponding {\it denominator} operation. A theorem of Schubert (See \cite{LKS}) gives a classification of 
the {\it rational knots} obtained as numerators of rational tangles. We do not need the full Schubert Theorem here but refer to \cite{KL2} for the reader who wants more information about it.\\

In Figure~\ref{UKK2} we illustrate the formula $N(T+ [n] + S + [-n]) = N(T + S).$ Actually this is simply an identity for tangles, but it is often occurring inside a numerator closure. Note that algebraically this is just the cancellation of $[n]$ with $[-n]$ but topologically it can be interpreted as the result of a rotation that involves the intermediate tangle $S.$ Such rotations occur in 
manipulating rope models and in the motions of curves in the force fields we have considered in the first section of the paper.\\

Now we are in position to look at Figure~\ref{UKK3} again. Note that the diagram indicated is $N((1,2,3) - (1,2)) = N([1] + (0,2,3) + [-1] + (0,-2))$ and so the $[1]$ and $[-1]$ cancel is in the previous paragraph. Again one can see this by an internal rotation in the numerator closure. The rest of the figure illustrates how this configuration is unknotted.\\

We can note that if $T = [a_1,a_2, \cdots , a_n]$ and $T^{!} = [a_n, a_{n-1}, \cdots , a_{1}],$ the tangle obtained by reversing all the terms in $T,$ then $N(T)$ and $N(T^{!})$ are topologically equivalent knots. The Figure~\ref{UKK4} will give the reader an appreciation of this fact. In that figure we illustrate how $N(2,3,4) = N(4,3,2)$ by representing each of these tangles in a 
partially braided form. The reader can take as an exercise to see that any continued fraction template for a rational tangle can be put into such a braided form.\\

Figure~\ref{UKK5} illustrates $N(T+S)$ for $T=(2,3,4)$ and $S=(3,2,2).$ We have the basic fact that the numerator of the sum of two rational tangles $T$ and $S$ is topologically equivalent to the numerator of a rational tangle $T \square S$ that is constructed from the two of them. Given $T = [a_1,a_2, \cdots , a_n]$ and $S = [b_1,b_2, \cdots ,b_{n-1},  b_n],$ let 
$T \square S$ be defined by the formula $$T \square S = (a_n, a_{n-1}, \cdots, a_2, a_1 + b_1, b_2, \cdots , b_{n-1}, b_n).$$ Then $N(T+S)$ is  topologically equivalent to $N(T \square S).$
In order to see this topological equivalence the Figure shows how to rotate and isotope $N(T + S)$ so that it is clearly equivalent to $N( (T \square S)^{!}).$ But we know that 
$N(T^{!}) = N(T)$ for any $T.$ So this gives the desired result. Thus in Figure~\ref{UKK5} we see that $$N((2,3,4) + (3,2,2)) = N(4,3,2+3,2,2) = N(4,3,5,2,2).$$\\

The rotation that we see in $N(T+S)$ in Figure~\ref{UKK5} is the key to producing many diagrams of the unknot that can be difficult for the dynamics of the program to undo.
The basic result is this: Let $T = [a_1,a_2, \cdots , a_n]$ and $S = [a_1,a_2, \cdots , a_{n-1}].$ Then $N(T - S)$ is topologically equivalent to the unknot. Figure~\ref{UKK3}
illustrates how this works. The tangle $S$ has just one less term than the tangle $T$ and is obtained from it by removing the last term of $T.$ We call $S$ the {\it truncate} of $T.$
In the figure the reader should be able to see that the $a_1$ from $T$ and the $-a_1$ from $-S$ flank a tangle and allow a rotation that cancels these terms and then brings $a_2$ and
$-a_2$ into cancelling position. This continues recursively and allows the knot to unwind and become unknotted. To see this result in tangle terms the following identity is of use:
$$(\cdots a,b,0,c,d,\cdots) = (\cdots a,b+c, d, \cdots)$$ where this equality can be read as topological equivalence of tangles, or as identity of continued fractions.
For example, in fractions we have 
$$(a,b,0,c,d) = a + 1/(b + 1/(0 + 1/(c + 1/d))))$$
$$ = a + 1/(b + 1/(1/(c + 1/d)))) $$
$$=  a + 1/((b + c) + 1/d)$$
$$ = (a,b+c, d).$$
The unknotting is then mirrored in a bit of algebra.
For example
$$N((a,b,c,d) - (a,b,c))  $$
$$= N(d,c,b,a-a,-b,-c)$$
$$ = N(d,c,b,0,-b,-c)$$
$$ = N(d,c,b-b,-c)$$
$$= N(d,c,0,-c)$$ 
$$= N(d,c-c)$$
$$ = N(d,0)$$
and $N(d,0)$ is unknotted.\\

By the same token, we can arrange tangle diagrams that collapse to certain knots.
For example,
$$N((1,1,1,c,b,a) -(1,1,1)) $$
$$= N(a,b,c,1,1,1-1,-1,-1)$$
$$= N(a,b,c,1,1,0,-1,-1)$$
$$ = N(a,b,c,0)$$
$$= N(a,b)$$
This means that we can give the dynamics program tangled versions of knots and unknots and perform experiments to see if the program 
will undo the unknots and collapse the complexified knots to simpler forms.\\

\subsection {Collapse Examples}

\noindent {\bf Notation.} We make the following definitions. Let $[[n]] = (1,1,1,\cdots, 1)$ with $n$ $1$'s in the continued fraction sequence.
Let $[n.m] = N([[n]] - [[m]]).$ This is the notation we have used in the first section of the paper. Thus Figure~\ref{Unknot} and Figure~\ref{Unknot1} give images
for $[11.10],$ a configuration that we know is unknotted from the above construction. These figures show that the KnotPlot program takes $[11.10]$ to a stable state
that is still entangled. The second instance is Figure~\ref{KnotMin} where we show how $[15.10]$ also has a complex stable state. By our discussion below we will see
that for $n > m$, $[n.m] = N([[n-m - 1]]).$ Thus $[15.10] = N([[6]])$ and this is a Figure Eight knot.The upshot of these constructions is that we can give
infinitely many examples of configurations that are either knotted or unknotted but have complex stable states in the force field algorithm of KnotPlot.\\

To prove that $[n.m] = N([[n-m - 1]]),$ consider an example:
$$[7.3] = N ((1,1,1,1,1,1,1)+(-1,-1,-1)) $$
$$= N(1,1,1,1,1,1,1-1,-1,-1) $$
$$= N(1,1,1,1,0) = N(1,1,1)$$
$$ = N([[7-3-1]]).$$
Thus the proof of this fact follows from the work we have already done, plus the use of this compressed notation.\\

We will also write $([[n]], a,b,c) = (1,1,1,\cdots 1,a,b,c)$ so that the $[[n]]$ can indicate $n$ consecutive $1$'s in the tangle.
Then we have the basic result that $$K_{n}(a,b,c,d) = N([[n+1]], a,b,c,d)-([[n]]))$$ is isotopic to $$N(a,b,c,d)$$  where $(a,b,c,d)$ can be replaced by any continued fraction tangle.
This result follows by the same reasoning as we have outlined above, and it provides a large collection of starting configurations for testing a self-repulsion program.\\

Our last example illustrates all the principles in this section. The basic knot we work with is the rational knot
$N(7,7,7,7)$ illustrated in its initial form in Figure~\ref{UKK7} and it its final fully repelled form in Figure~\ref{UKK8}.
In Figure~\ref{UKK9} we show the initial form of $K_{11}(7,7,7,7) = N(([[12]],7,7,7,7) - ([[11]]).)$
In Figure~\ref{UKK10} we show the final self-repelled form that occurs from this beginning.
We see that in this case $K_{11}(7,7,7,7)$ does not reduce to the minimal form of $N(7,7,7,7)$ in the KnotPlot environment.
It is our conjecture that no environment of this type will fully reduce all of the configurations $K_{n}(7,7,7,7),$ and we challenge users or designers of other 
programs to experiment with this family of configurations of the knot $N(7,7,7,7).$ There is, of course, nothing special about the choices we have made here.
The phenomenon of non-collapse will occur for sufficiently large complexifications of any knot or unknot.\\

\end{document}